\documentclass[a4paper,12pt,reqno]{article}

\usepackage[T1]{fontenc}
\usepackage{amsmath}
\usepackage{amsthm}
\usepackage[author-year,nobysame]{amsrefs}
\usepackage{amsfonts}
\usepackage{mathrsfs}
\usepackage{exscale}
\usepackage{scrtime}

\numberwithin{equation}{section}

\newcommand{\set}[1]{\{#1\}}

\theoremstyle{plain}

\newtheorem{theorem}{Theorem}[section]

\theoremstyle{definition}

\theoremstyle{remark}

\newcommand{\prob}{\mathbb P}
\newcommand{\pp}[1]{{\mathbb P}[#1]}

\newcommand{\law}[1]{{\mathscr L}(#1)}
\newcommand{\biglaw}[1]{{\mathscr L}\bigl(#1\bigr)}

\newcommand{\bigcondlaw}[2]{{\mathscr L}\bigl(#1\bigm|#2\bigr)}

\newcommand{\ee}{{\mathbb E}}

\newcommand{\dtv}[2]{d_{\scriptscriptstyle {\mathrm TV}}(#1,#2)}
\newcommand{\bigdtv}[2]{d_{\scriptscriptstyle {\mathrm TV}}\bigl(#1,#2\bigr)}

\DeclareMathOperator{\be}{Be}

\DeclareMathOperator{\po}{Po}

\newcommand{\ddens}{p_{\theta}}

\newcommand{\nn}{{\mathbb N}}
\newcommand{\zz}{{\mathbb Z}}
\def\integ{\zz}
\newcommand{\zzp}{{\mathbb Z}_+}

\newcommand{\lo}[1]{{\mathrm o}(#1)}

\newcommand{\bigbo}[1]{{\mathrm O}\bigl(#1\bigr)}

\newcommand{\eqa}{\begin{eqnarray}}
\newcommand{\ena}{\end{eqnarray}}
\newcommand{\eqs}{\begin{eqnarray*}}
\newcommand{\ens}{\end{eqnarray*}}
\newcommand{\eq}{\begin{equation}}
\newcommand{\en}{\end{equation}}

\def\Ref#1{(\ref{#1})}
\def\Le{\ \le\ }
\def\Def{\ :=\ }

\def\Bl{\left(}
\def\Br{\right)}
\def\Blb{\left\{}
\def\Brb{\right\}}

\def\ep{\hfill$\Box$\bsk}

\def\a{\alpha}

\def\th{\theta}

\def\r{\rho}

\def\sn{\sum_{i=1}^n}

\def\Eq{\ =\ }
\def\pr{\prob}
\def\Po{\po}
\def\nat{\nn}

\def\un{^{(n)}}

\def\ignore#1{}

\def\half{{\textstyle \frac12}}

\def\ep{\hfill$\Box$}

\def\bP{{\overline P}}
\def\sjan{\sum_{j=a+1}^n}
\def\giv{\,|\,}
\def\tV{{\widetilde V}}
\def\hZ{{\widehat Z}}
\def\hT{{\widehat T}}

\newcommand{\cn}{C^{\scriptscriptstyle (n)}} 

\begin{document}

\title
{Couplings for irregular combinatorial assemblies}

\author{A. D. Barbour\footnote{Angewandte Mathematik, Universit\"at Z\"urich,  
Winterthurertrasse 190, CH-8057 Z\"URICH;  E-mail {\tt a.d.barbour@math.uzh.ch};}
and A. P\'osfai\footnote{Department of Mathematics, Tufts University, 503 Boston Avenue, Medford, MA 02155, USA
and Analysis and Stochastics Research Group of the Hungarian Academy of Sciences, Bolyai Institute, 
University of Szeged, Aradi v\'ertan\'uk tere 1, Szeged 6720, Hungary; E-mail {\tt anna.posfai@tufts.edu}}\\
Universit\"at Z\"urich and Tufts University}

\date{}
\maketitle

\begin{abstract}
When approximating the joint distribution of the component counts of a 
decomposable combinatorial structure that is `almost' in the logarithmic
class, but nonetheless has irregular structure, it is useful to be able
first to establish that the distribution of a certain sum of non-negative integer 
valued random variables is smooth.  This distribution is not like the
normal, and individual summands can contribute a non-trivial amount to the
whole, so its smoothness is somewhat surprising.  In this paper, we 
consider two coupling approaches to establishing the smoothness, and
contrast the results that are obtained.
\end{abstract}

{\it Running head:}  Couplings for combinatorial assemblies

\section{Introduction}

Many of the classical random decomposable combinatorial structures, such
as random permutations and random polynomials over a finite field, 
have component structure satisfying
a \emph{conditioning relation}: if~$\cn_i$ denotes the number of 
components of size~$i$, the distribution of the vector of component counts
$(\cn_1, \dotsc, \cn_n)$ of a structure of size~$n$ can be expressed as
\begin{equation}\label{AB-cond-rel}
    \biglaw{\cn_1, \dotsc, \cn_n} = \bigcondlaw{Z_1, \dotsc ,Z_n}{T_{0,n}=n} \, ,
\end{equation}
where $(Z_i,\,i\ge1)$ is a fixed sequence of independent non-negative integer 
valued random variables, and $T_{a,n} := \sum_{i=a+1}^n iZ_i$, $0\le a < n$.  
If, as in the examples above, the~$Z_i$ also satisfy
\begin{equation}\label{AB-logarithmic}
  i\pp{Z_i=1} \ \to\ \theta  
     \qquad \text{and} \qquad \th_i \Def i\ee Z_i\ \to\ \th \, ,
\end{equation}
the combinatorial structure is called \emph{logarithmic}.  
It is shown in Arratia, Barbour \& Tavar\'{e}~(2003) [ABT] that
combinatorial structures satisfying the conditioning relation and
slight strengthenings of the logarithmic condition share many common
properties.  For instance, if~$L\un$ is the size of the largest component,
then 
\eq\label{big}
   n^{-1}L\un \to_d L,
\en
where~$L$ has probability density function
$f_{\theta}(x)  := e^{\gamma \theta} \Gamma(\theta+1) x^{\theta-2} \ddens((1-x)/x)$, 
$x \in (0,1]$, and~$p_\th$ is the density of the Dickman 
distribution~$P_\th$ with parameter~$\th$, given in Vervaat~(1972, p.~90).
Furthermore, for any sequence $(a_n,\,n\ge1)$ with $a_n = o(n)$, 
\begin{equation}\label{small}
   \lim_{n \to \infty} 
      \bigdtv{\law{\cn_1,\dotsc,\cn_{a_n}}}{\law{Z_1,\dotsc,Z_{a_n}}} = 0 \, .
\end{equation}
Both of these convergence results can be complemented by estimates of the
approximation error, under appropriate conditions.

If the logarithmic condition is not satisfied, as in certain of the
additive arithmetic semigroups introduced in Knopfmacher~(1979), the
results in [ABT] are not directly applicable.  However, in Manstavi\v cius~(2009)
and in Barbour \& Nietlispach~(2010) [BN], it is shown that the logarithmic
condition can be relaxed to a certain extent, without disturbing the validity
of~\Ref{small}, and that~\Ref{big} can also be recovered, if the convergence
in~\Ref{AB-logarithmic} is replaced by a weaker form of convergence.  
A key step in the proofs of these results is to be able to show
that, for sequences $a_n = \lo n$, the normalized sum~$n^{-1}T_{a_n,n}$ converges both
in distribution and locally to the Dickman distribution~$P_\th$, 
and that the error rates in these approximations can be controlled. To do so,
it is in turn necessary to be able to show that, under suitable conditions, 
\begin{equation}\label{dtv to zero intro}
  \lim_{n \to \infty} \bigdtv{\law{T_{a_n,n}}}{\law{T_{a_n,n} + 1}} \Eq 0 \, , 
      \qquad \text{for all $a_n=\lo{n}$,}
\end{equation} 
and that the error rate can be bounded by a power of~$\{(a_n+1)/n\}$. In this note,
we explore ways of using coupling to prove such estimates, in the
simplest case in which the~$Z_i$ have Poisson distributions.  The first
of these, an improvement over the Mineka coupling, was introduced in~[BN].
It is extremely flexible in obtaining error rates bounded by a power 
of~$\{(a_n+1)/n\}$ for a wide variety of choices of the means~$\th_i$,
and it is in no way restricted to Poisson distributed~$Z_i$'s.
Here, we show that, despite its attractions, it does not achieve the
best possible error rate under ideal circumstances.  The second approach
works only in much more restricted situations, but is then capable of attaining
the theoretically best results. 

In the case of Poisson distributed~$Z_i$, the distribution of~$T_{an}$ is
a particular compound Poisson distribution, with parameters determined
by~$n$ and by the~$\th_i$, and it is tempting to try to
approximate the distribution of~$n^{-1}T_{a_n,n}$ by first approximating
by the distribution that would be obtained if $\th_i = \th$ for all~$i$.
A natural way of obtaining compound Poisson approximation is then to use
Stein's method (Barbour, Chen \& Loh~1992). Difficulties arise, however,
because the conditions of their Theorem~5 (needed to get useful bounds on
the solution to the Stein equation) are not satisfied unless $a=0$,
and, even then, the bounds obtained are not as useful as they might be;
better information for this particular case can be found in [ABT, Chapter~9].  
And, even using this approach, it still seems necessary first to bound the error
in~\Ref{dtv to zero intro}, in order to obtain useful results.

\section{A Mineka--like coupling}\label{coupling subsection}

Let $\set{X_i}_{i \in \nn}$ be mutually independent $\zz$-valued random variables, and let
$ S_n := \sn X_i$. 
The Mineka coupling, developed independently by Mineka~(1973) and R\"osler~(1977)
(see also Lindvall~(2002, Section II.14)) yields a bound of the form
\begin{equation}\label{barbour xia bound}
  \bigdtv{\law{S_n}}{\law{S_n + 1}} \Le \Bigl(\frac\pi2 \sum\nolimits_{i=1}^{n} u_i \Bigr)^{-1/2} \, ,
\end{equation}
where
\begin{equation*}
   u_i \Def  \Bigl( 1 - \bigdtv{\law{X_i}}{\law{X_i+1}}\Bigr) \, ;
\end{equation*}
see Mattner \& Roos~(2007, Corollary~1.6). The proof is based on coupling copies
 $\{X_i'\}_{i \in \nn}$ and $\{X_i''\}_{i \in \nn}$ of $\{X_i\}_{i \in \nn}$ 
in such a way that
\begin{equation*}
  V_n \Def \sum_{i=1}^n \bigl(X_i - X_i'\bigr)\,, \qquad n \in \nn,
\end{equation*}
is a symmetric random walk with steps in $\{-1,0,1\}$. Writing $S'_i := 1 + \sum_{j=1}^i X_j' \sim
S_i+1$
and $S_i'' := \sum_{j=1}^i X_j'' \sim S_i$, 
so that $V_i + 1 = S'_i - S_i''$,
the coupling inequality (Lindvall 1992, Section~I.2) then shows that
\begin{equation*}
   \bigdtv{\law{S_n}}{\law{S_n + 1}} \Le \pp{\tau > n} \Eq \pp{V_n = \{-1,0\}} \, ,
\end{equation*}
where $\tau$ is the time at which $\{V_n\}_{n \in \zzp}$ first hits level $-1$,
and the last equality follows from the reflection principle.
However, this inequality gives slow convergence rates, if $X_i=i Z_i$ and the $Z_i$ 
are as described in the Introduction;  typically, $\bigdtv{\law{iZ_i}}{\law{iZ_i+1}}$
is equal to~$1$,
and, if $X_i$ is taken instead to be $(2i-1)Z_{{2i-1}} + 2iZ_{2i}$, we still expect to have
$1 - \bigdtv{\law{X_i}}{\law{X_i+1}} \asymp i^{-1}$, leading to bounds of the form
\begin{equation}\label{Mineka-error}
   \bigdtv{\law{T_{a_n,n}}}{\law{T_{a_n,n}+1}} = \bigbo{(\log (n/\{a_n+1\}))^{-1/2}} \, .
\end{equation}

The reason that the Mineka coupling does not work efficiently in our setting is
that,  once the random walk~$V_n$ takes some value~$k$, it has to achieve a
preponderance of~$k+1$ negative steps, in order to get to the state~$-1$, and this
typically requires many jumps to realize.  Since, at the $i$-th step, the probability of there
being a jump is of order~$i^{-1}$, it thus takes a very long time for such an event to
occur, and the probability of this not happening before time~$n$ is then relatively
large. In [BN], the difficulty is overcome by observing that the Mineka random
walk can be replaced by another Markov chain~$(\tV_n, n\ge1)$, still constructed 
from copies $(Z'_i,\,i\ge1)$ and~$(Z_i'',\,i\ge1)$ of the original sequence, but now 
associated differently with one another.  The basic idea is to note that, if~$\tV_i=k$,
then the random variables $X_{i+1}' := jZ_j' + (j+k+1)Z_{j+k+1}'$ and
$X_{i+1}'' := jZ_j'' + (j+k+1)Z_{j+k+1}''$ can be coupled in such a way that
$X_{i+1}' - X_{i+1}'' \in \{-(k+1),0,(k+1)\}$, for any~$j$ such that the indices $j$ and~$j+k+1$
have not previously been used in the construction.  Hence a single jump has probability
$1/2$ of making~$\tV$ reach~$-1$. The construction starts as for 
the Mineka walk, but if the first jump takes~$\tV$ to~$+1$, then the chain switches
to jumps in $\{-2,0,2\}$; and subsequently, if $\tV$ is in the state $k=2^r-1$, the chain
makes jumps in $\{-2^r,0,2^r\}$.  Clearly, this construction can be used with
$Z_i \sim \Po(i^{-1}\th_i)$, even when many of the~$\th_i$ are zero.  A number of settings
of this kind are explored in detail in [BN]; for instance, when $\th_i \ge \th^*$
for all~$i$ in $\{r\integ_++t\}\cup\{s\integ_++u\}$, where $r$ and~$s$ are coprime.  
Very roughly, provided that a non-vanishing
fraction of the~$\th_i$ exceed some fixed value~$\th_* > 0$, the probability that
$\tV$ reaches~$-1$ before time~$n$ is of order~$n^{-\a}$, for some $\a > 0$, an
error probability exponentially smaller than that in~\Ref{Mineka-error}.

Here, we make the following observation.  Suppose that we have the
ideal situation in which $\th_i = \th^* > 0$
for {\em every}~$i$.   Then the probability that a coupling, constructed
as above, should fail is at least of magnitude~$n^{-\th^*/2}$.  In 
Section~\ref{Poisson}, it is shown that the total variation distance 
in~\Ref{dtv to zero intro}  is actually of order $n^{-\min\{\th^*,1\}}$ under these
circumstances, so that the estimates of this distance obtained
by the~[BN] coupling are typically rather weaker.  It is thus of interest
to find ways of attaining sharper results.  The coupling given in Section~\ref{Poisson}
is one such, but it is much less widely applicable.

The coupling approach given in~[BN] evolves by 
choosing a pair of indices $M_{i1}<M_{i2}$ at each step~$i$, with the choice
depending on the values previously used: no index can be used more than once, and 
$M_{i2}-M_{i1} = \tV_{i-1}+1$, so that one jump in the right direction
leads immediately to a successful coupling. Then, if $(M_{i1},M_{i2}) = (j,j+k+1)$,
the pair $X_i'$ and~$X_i''$ is constructed as above, by way of copies of the random variables
$jZ_j$ and $(j+k+1)Z_{j+k+1}$.  The probability of a jump taking place is then
roughly $2\th^*/(j+k+1)$, and, if a jump occurs, it has probability~$1/2$ of
taking the value $-(k+1)$, leading to success.   The main result of this
section is the following lower bound for the failure probability of such a procedure.

\begin{theorem}\label{BN-coupling}
For any coupling constructed as above, the 
probability~$\pr[F]$ that the coupling is not successful is bounded
below by
$$
   P[F] \ \ge\ \prod_{i=1}^{\lfloor n/2 \rfloor} (1 - \half\min\{\th^*/i,1/e\})
     \ \asymp\ n^{-\th^*/2}.
$$
\end{theorem}

\proof
In order to prove the lower bound, we couple two processes, one of which makes
more jumps than the other. We start by letting $(U_i,\,i\ge1)$ be independent uniform random
variables on $[0,1]$. The first process is much as discussed above.  It is
defined by a sequence of pairs of indices $M_{i1}<M_{i2}$, $1\le i\le I^*$, from~$[n]
:= \{i \in \nat\colon\,i \le n\}$,
with $I^* \le \lfloor n/2 \rfloor$ the last index for which a suitable pair
can be found. No index is ever used twice, and the choice of $(M_{i1},M_{i2})$ 
is allowed to depend on $((M_{j1},M_{j2}, U_j),\,1\le j < i)$. 
We set $Y_i = I[U_i \le p(M_{i1},M_{i2})]$, where
\[
    p(m_1,m_2) \Def 2e^{-\th^*/m_1}\,(\th^*/m_2) e^{-\th^*/m_2} \ <\ 1/e,
\]
for $m_1 < m_2$, representing the indicator of a jump of~$\pm (m_2-m_1)$ being made 
by the first process at time~$i$.  For the second, we inductively define 
$R_{i} := \{\r(1),\ldots,\r(i)\}$ by taking $R_0 = \emptyset$ and
\[
    \r(i) \Def \max\{r \in [n/2]\setminus R_{i-1}\colon\,2r \le M_{i2}\};
\]
we shall check at the end of the proof that~$\r(i)$ always exists. 
(The second process, that we do not really need in detail, uses the pair 
$(2\r(i)-1,2\r(i))$ at stage~$i$.)
We then define $Z_i := I[U_i \le \min\{\th^*/\r(i),1/e\}]$, noting that
$p(M_{i1},M_{i2}) \le \min\{\th^*/\r(i),1/e\}$, entailing $Z_i \ge Y_i$ a.s.\ for
all~$i$.  Finally, let $(J_i,\,i\ge1)$ be distributed as~$\be(1/2)$, independently
of each other and everything else.

The event that the first process makes no successful jumps can be described as the
event
\[
    F \Def \Blb \sum_{i=1}^{I^*} Y_iJ_i = 0 \Brb.
\]
We thus clearly have
\[
    F \ \supset\ \Blb \sum_{i=1}^{\lfloor n/2 \rfloor} Z_i J_i = 0 \Brb,
\]
where, for $I^* < i \le n/2$, we take $\r(i) := \min\{r\in [n/2]\setminus R_{i-1}\}$,
and $R_i := R_{i-1}\cup\{\r(i)\}$.  But now the~$Z_i$, suitably reordered, are
just independent Bernoulli random variables with means $\min\{\th^*/r,1/e\}$,
$1\le r\le n/2$, and hence
\[
   P[F] \ \ge\ \prod_{i=1}^{\lfloor n/2 \rfloor}  (1 - \half\min\{\th^*/i,1/e\})
     \ \asymp\ n^{-\th^*/2}.
\]

It remains to show that the~$\r(i)$ are well defined at each stage, which requires
that 
\[
  S_i \Def \{r\in[n/2]\setminus R_{i-1}\colon\,2r \le M_{i2}\}\ \neq\ \emptyset,
\]
$1\le i\le I^*$.  For $i=1$, $m_{12} \ge 2$, so the start is successful.  Now,
for $2\le i\le n/2$, suppose that
\[
    r(i-1) \Def \max\{s\colon\, R_{i-1} \supset \{1,2,\ldots,s\}\}.
\]
Then $1,2,\ldots,r(i-1)$ can be expressed as
$\r(i_1),\r(i_2),\ldots,\r(i_{r(i-1)})$, for some indices $i_1,i_2,\ldots,i_{r(i-1)}$.
For these indices, we have $M_{i_l,2} \le 2r(i-1)+1$, $1\le l\le r(i-1)$, 
since $r(i-1)+1 \notin R_{i-1}$ and, from the definition of $\r(\cdot)$, 
we could thus not choose $\r(i_l) \le r(i-1)$ if $M_{i_l,2} \ge 2r(i-1)+2$.
Hence, also, $M_{i_l,1} \le 2r(i-1)+1$, and, because all the $M_{is}$ are
distinct, $\{M_{i_l,s},\, 1\le s\le 2, 1 \le l \le r(i-1)\}$ is a set
of $2r(i-1)$ elements of $[2r(i-1)+1]$. Thus,
when choosing the pair $(M_{i1},M_{i2})$, there is only at most one element
of $[2r(i-1)+1]$ still available for choice, from which it
follows that $M_{i2} \ge 2r(i-1) + 2$: so $r(i-1)+1 \in S_i$, and hence~$S_i$
is not empty.
\ep

\section{A Poisson--based coupling}
\label{Poisson}
In this section, we show that a coupling can be constructed that gives good
error rates in~\Ref{dtv to zero intro} when  $Z_j \sim \Po(j^{-1}\th^*)$, for some fixed
$\th^* > 0$.  If $Z_j \sim \Po(j^{-1}\th_j)$
with $\th_j \ge \th^*$,  the same order of error can immediately be deduced 
(though it may no longer be optimal), since, for Poisson
random variables, we can write $T_{an} = T^*_{an} + T'$, with $T^*_{an}$
constructed from independent random variables~$Z^*_j \sim \Po(j^{-1}\th^*)$, 
and with~$T'$ independent of~$T^*_{an}$.  

Because of the
Poisson assumption, the distribution of $T_{an} := \sjan jZ_j$ can 
equivalently be
re-expressed as that of a sum of a random number $N \sim \Po(\th^*h_{an})$
of independent copies of a random variable~$X$ having $\pr[X=j] = 1/\{jh_{an}\}$,
$a+1 \le j \le n$, where $h_{an} := \sjan j^{-1}$.  Fix $c > 1$, define $j_r :=
\lfloor c^r \rfloor$, and set 
\[
    r_0 \Def r_0(a) \Def \lceil \log_c(a+1) \rceil,\qquad 
        r_1 \Def r_1(n) \Def \lfloor \log_c n \rfloor.
\]
Define independent random variables $(X_{ri},\,r_0 \le r < r_1,\, i\ge1)$ and 
$(N_r,\,r_0 \le r < r_1)$, with $N_r \sim \Po(\th^*h_{an}p_r)$ and
\[
    \pr[X_{ri} = j] \Eq 1/\{jh_{an}p_r\},\qquad j_r \le j < j_{r+1},
\]
where
\[
    p_r \Def \sum_{j=j_r}^{j_{r+1}-1} \frac1{jh_{an}};
\]
define $\bP_r := \sum_{s=r}^{r_1-1} p_s \le 1$.
Then we can write~$T_{an}$ in the form
\[
    T_{an} \Eq Y + \sum_{r=r_0}^{r_1-1} \sum_{i=1}^{N_r} X_{ri},
\]
where $Y$ is independent of the sum; the~$X_{ri}$ represent the realizations
of the copies of~$X$ that fall in the interval $C_r := [j_r,j_{r+1})$, and~$Y$ accounts 
for all $X$-values not belonging to one of these intervals.  The idea is then to
construct copies $T'_{an}$ and~$T''_{an}$ of~$T_{an}$ with $T'_{an}$ coupled to~$T''_{an}+1$,
by using the same~$N_r$ for both, and trying to couple
one pair $X'_{ri}$ and $X''_{ri}+1$ exactly, declaring
failure if this doesn't work.  Clearly, such a coupling can only be attempted
for an~$r$ for which $N_r \ge 1$.  Then exact coupling can be achieved between
$X'_{r1}$ and $X''_{r1}+1$ with probability $1 - 1/\{j_rh_{an}p_r\}$, since the
point probabilities for~$X_{r1}$ are decreasing. Noting that the $p_r$ are all of the same magnitude, 
it is thus advantageous to try
to couple with~$r$ as large as possible.  This strategy leads to the following theorem.

\begin{theorem}\label{Poisson-coupling}
   With $Z_j \sim \Po(j^{-1}\th^*)$, $j\ge1$, we have
\[ 
    \bigdtv{\law{T_{an}}}{\law{T_{an} + 1}} \Eq O(\{(a+1)/n\}^{\th^*} + n^{-1}),
\]
if $\th^*\neq1$; for $\th^*=1$,
\[ 
    \bigdtv{\law{T_{an}}}{\law{T_{an} + 1}} \Eq O(\{(a+1)/n\} + n^{-1}\log\{n/(a+1)\}).
\]
\end{theorem}

\proof
We begin by defining 
\[
   B_r \Def \Bl \bigcap_{s=r+1}^{r_1-1}\{N_s=0\} \Br \cap \{N_r\ge1\},\qquad r_0 \le r < r_1,
\]
and setting $B_0 := \bigcap_{s=r_0}^{r_1-1}\{N_s=0\}$. 
On the event~$B_r$, write
$X''_{r1} = X'_{r1}-1$ if $X'_{r1} \ne j_r$, with $X''_{r1}$ so distributed on
the event $A_r := \{X'_{r1} = j_r\}$ that its overall distribution is correct.
All other pairs of random variables $X'_{r'i}$ and~$X''_{r'i}$, 
$(r',i)\in ([r_0,\ldots,r_1-1]\times\nat)\setminus\{(r,1)\}$, are set to be
equal on~$B_r$.  This generates copies $T'_{an}$ and $T''_{an}$ of~$T_{an}$,
with the property that $T'_{an} = T''_{an}+1$, except on the event
\[
     E \Def B_0 \cup \Bl \bigcup_{r=r_0}^{r_1-1} (B_r \cap A_r) \Br.
\]

It is immediate from the construction that 
\[
   \pr[B_r] \Eq \exp\{-\th^* h_{an} \bP_{r+1}\}(1 - e^{-\th^* h_{an} p_r}),
     \quad r_0 \le r < r_1,
\]
and that $\pr[B_0] = \exp\{-\th^* h_{an} \bP_{r_0}\}$; and
$\pr[A_r \giv B_r] = 1/\{j_rh_{an}p_r\}$.  This gives all the ingredients
necessary to evaluate the probability
\[
   \pr[E] \Eq \pr[B_0] + \sum_{r=r_0}^{r_1-1} \pr[B_r] \pr[A_r \giv B_r].
\]
In particular, as $r\to\infty$, $j_r \sim c^r$, $h_{an}p_r \sim \log c$ 
and $h_{an}\bP_{r+1} \sim (r_1(n)-r)\log c$, from which it follows that
$\pr[B_r] \sim c^{-\th^*(r_1(n)-r)}(1-c^{-\th^*})$, $\pr[A_r \giv B_r] \sim 
1/\{c^r\log c\}$ and
\[
   \pr[B_0] \ \asymp\ c^{-\th^*(r_1(n) - r_0(a))} \ \asymp\ \{(a+1)/n\}^{\th^*}.
\]
Combining this information, we arrive at
\[
   \pr[E] \ \asymp\ \{(a+1)/n\}^{\th^*} + \sum_{r=r_0(a)}^{r_1(n)-1} c^{-r}\,c^{-\th^*(r_1(n)-r)}.
\]
For $\th^* > 1$, the dominant term in the sum is that with $r = r_1(n)-1$, and
it follows from the definition of~$r_1(n)$ that then 
\[
   \pr[E] \ \asymp\ \{(a+1)/n\}^{\th^*} + c^{-r_1(n)} \asymp \{(a+1)/n\}^{\th^*} + n^{-1}.
\]
For $\th^* < 1$, the dominant term is that with $r=r_0(a)$, giving
\[
   \pr[E] \ \asymp\ \{(a+1)/n\}^{\th^*} + n^{-\th^*}(a+1)^{-(1-\th^*)} \ \asymp\ \{(a+1)/n\}^{\th^*}.
\]
For $\th^*=1$, all terms in the sum are of the same order, and we get
$$
  {}\hskip1.52in\pr[E]\ \asymp\ \{(a+1)/n\} + n^{-1}\log(n/(a+1)).\hskip1.52in\Box
$$

\bigskip
Note that the element $\{(a+1)/n\}^{\th^*}$ appearing in the errors is very easy to
interpret, and arises from the probability of the event that $T_{an} = 0$, a
value unattainable by $T_{an}+1$.  Furthermore, the random variable~$T_{an}$
has some point probabilities of magnitude~$n^{-1}$ [ABT, p.~91], so that~$n^{-1}$
is always a lower bound for the order of $\bigdtv{\law{T_{an}}}{\law{T_{an} + 1}}$. 
Hence the order of approximation in Theorem~\ref{Poisson-coupling} is 
best possible if $\th^*\ne1$.  However, for $a=0$ and $\th^*=1$, the point probabilities 
of $T_{0n}$ are decreasing, and since their maximum is of order~$O(n^{-1})$, 
the logarithmic factor in the case~$\th^*=1$ is not sharp, at least for $a=0$.

The method of coupling used in this section can be extended in a number of ways.
For instance, it can be used for random variables~$Z_j$
with distributions other than Poisson, giving the same order of error as long as 
$\dtv{\law{Z_j}}{\Po(\th_j/j)} = O(j^{-2})$.  This is because, first,
for some $K < \infty$,
\[
   \pr[B_r] \Le K \exp\{-\th^* h_{an} \bP_{r+1}\}(1 - e^{-\th^* h_{an} p_r}),
     \quad r_0 \le r < r_1,
\]
and $\pr[B_0] \le K \exp\{-\th^* h_{an} \bP_{r_0}\}$, where, in the definitions
of the~$B_r$, the events $\{N_s=0\}$ are replaced by $\{Z_j=0,\,j_s \le j < j_{s+1}\}$.
Secondly, we immediately have
\[
   \dtv{\law{(Z_j,\,j_r \le j < j_{r+1})}}{(\hZ_j,\,j_r \le j < j_{r+1})} \Eq O(j_r^{-1}),
\]
where the $\hZ_j\sim \Po(\th_j/j)$ are independent, and hence that
\[   
   \dtv{\law{T_{j_r-1,j_{r+1}-1}}}{\law{\hT_{j_r-1,j_{r+1}-1}}} \Eq O(j_r^{-1}),
\]
where~$\hT_{rs}$ is defined as~$T_{rs}$, but using the~$\hZ_j$. Thus, on the
event~$B_r$, coupling can still be achieved except on an event of probability
of order~$O(j_r^{-1})$.

It is also possible to extend the argument to allow for gaps between the
intervals on which $\th_j \ge \th^*$.  Here, for $0 < c_1 \le c_2$, the
intervals $[j_r,j_{r+1}-1]$ can be replaced by intervals $[a_r,b_r]$,
such that $b_r/a_r \ge c_1$ and $a_r \ge kac_2^r$ for some~$k$ and for each $1 \le r \le R$,
say.  The argument above then leads to a failure probability of at most
\[
   O\Bl c_1^{-R\th^*} + \sum_{r=1}^{R} \frac1{ac_2^{r}}\,c_1^{-\th^*(R-r)} \Br.
\]
If $c_1^{\th^*} > c_2$, the failure probability is thus at most of order 
$O(c_1^{-R\th^*} + 1/\{ac_2^R\})$; if $c_1^{\th^*} < c_2$, it is of
order $O(c_1^{-R\th^*})$.  In Theorem~\ref{Poisson-coupling} above, we
have $c_1=c_2=c$, $k=1$ and $c^{-R} \asymp (a+1)/n$, and the results are
equivalent.

However, the method is still only useful if there are
long stretches of indices~$j$ with~$\th_j$ uniformly bounded below.  This
is in contrast to that discussed in the previous section, which is flexible
enough to allow sequences~$\th_j$ with many gaps.  It would be interesting to
know of other methods that could improve the error bounds obtained by these
methods.

\end{document}